\documentclass[twoside,11pt]{article}

\usepackage{jmlr2e}
\usepackage{amsfonts}
\usepackage{amsmath, float, subcaption, booktabs, physics}
\usepackage{amssymb,amsmath,wrapfig}
\usepackage{cancel}

\usepackage[all]{xy}
\usepackage{url}
\usepackage{listings}
\usepackage{color}

\renewcommand{\tilde}{\widetilde}

\usepackage{tikz}
\usepackage{dsfont}
\usepackage[inline]{enumitem}
\usepackage{tikz-3dplot}
\usetikzlibrary{shapes.geometric, calc}
\graphicspath{ {./final_project_images/} }
\newcommand{\comments}[1]{}     
\usepackage{hyperref}           
\usepackage{algorithm}
\usepackage[noend]{algpseudocode}

\newcommand{\R}{\mathbb{R}}

\usepackage{amsbsy,amscd,amsgen}

\usepackage{appendix}

\usepackage{comment}
\usepackage{multirow}
\usepackage{bm}

\ShortHeadings{Convolutional Autoencoders for Reduced-Order Modeling}{Venkat, Smith, and Kelley}
\firstpageno{1}

\begin{document}

\title{Convolutional Autoencoders for Reduced-Order Modeling}

\author{\name Sreeram Venkat \email srvenkat@ncsu.edu \\
       \addr Department of Mathematics\\
       North Carolina State University\\
       Raleigh, NC 27695, USA
       \AND
       \name Ralph C. Smith \email rsmith@ncsu.edu \\
       \addr Department of Mathematics\\
       North Carolina State University\\
       Raleigh, NC 27695, USA
       \AND
       \name Carl T. Kelley \email tim\_kelley@ncsu.edu\\
       \addr Department of Mathematics\\
       North Carolina State University\\
       Raleigh, NC 27695, USA
       }

\editor{}

\maketitle

\begin{abstract}
In the construction of reduced-order models for dynamical systems, linear projection methods, such as proper orthogonal decompositions, are commonly employed. However, for many dynamical systems, the lower dimensional representation of the state space can most accurately be described by a \textit{nonlinear} manifold. Previous research has shown that deep learning can provide an efficient method for performing nonlinear dimension reduction, though they are dependent on the availability of training data and are often problem-specific \citep[see][]{carlberg_ca}. Here, we utilize randomized training data to create and train convolutional autoencoders that perform nonlinear dimension reduction for the wave and Kuramoto-Shivasinsky equations. Moreover, we present training methods that are independent of full-order model samples and use the manifold least-squares Petrov-Galerkin projection method to define a reduced-order model for the heat, wave, and Kuramoto-Shivasinsky equations using the same autoencoder.

\end{abstract}

\begin{keywords}
  Autoencoders, Reduced-Order Modeling, Nonlinear Projection, Manifold Learning, Least-Squares-Petrov-Galerkin
\end{keywords}

\section{Introduction}
The full-scale models, which are often used to model physical phenomena, are generally systems of PDEs derived from natural laws. However, we cannot always use the full-scale models in practice as they are usually very computationally intensive. Specifically, applications such as Bayesian inference for parameter estimation, uncertainty propagation, experimental design, and real-time control require thousands to millions of model evaluations. Computational costs are compounded when analyzing large-scale systems through coupled models. To remedy this issue, we create reduced-order models, or surrogate models, that reduce computational intensity by using techniques including interpolation, projection, and Gaussian processes \citep[see][]{interpolatory, projection, gaussian-processes-ML, gaussian-processes, high-dim-model-reduction}. As data-driven analysis of dynamical systems becomes more prevalent, a common solution has been to construct a reduced order model (ROM) by projection onto a linear subspace of the original state space using methods such as proper orthogonal decomposition. In many cases, however, the properties of the original dynamical system are most accurately quantified by a projection onto a \textit{nonlinear} manifold \citep[see][]{carlberg_ca}.

Determining a nonlinear manifold that best fits the given data does not have a simple analytic solution that is computationally practical for large-scale applications. Researchers such as \citet{carlberg_ca} have used \textbf{convolutional autoencoders} from deep learning for this purpose. Initially developed for image processing, convolutional autoencoders have been shown to provide a computationally efficient method for the nonlinear dimension reduction of dynamical systems. However, one of the issues that they face is the problem of obtaining training data. Unlike image processing autoencoders, which have a virtually limitless amount of training data, autoencoders applied to physical systems face a severe shortage of training data. The training data for the dimension reduction method primarily comes from sampling the full-order model; in practice, it may be that one can only obtain 10-100 samples of the full-order model. As a result, we need to devise a training method for these autoencoders that does not depend on having many samples of the full-order model. Moreover, to remain computationally efficient, we would like the training of the autoencoder to happen completely offline. In the subsequent discussion, we present methods for developing and training autoencoders that have no dependence on full-order model samples, can be trained fully offline, and can be applied to completely different problems.

\section{Mathematical Models}\label{sec:models}
In this paper, we consider three mathematical models of dynamical systems: the heat equation, the wave equation, and the Kuramoto-Shivasinsky equation. We outline here these models and their associated numerical solutions.

\subsection{Definitions}\label{subsec:defs}
 The heat equation is the classical example of a parabolic PDE used to model physical phenomena such as diffusion and heating/cooling.We consider the heat equation 
 
 \begin{gather}
    \begin{cases}
        u_{t} = u_{xx},& x\in \Omega, t> 0,\\
        u(0,x) = u_0(x) = 5x(1-x),& x\in \Omega,\\
        u(t,0) = u(t,1) = 0,& t>0\\
    \end{cases}\label{eq:heat}
\end{gather}

\noindent
 on the bounded domains on $\Omega = [-1,1]$ and $[0,1]$. Here, we choose the constants $\alpha_1,\alpha_2, \beta_1,$ and $\beta_2$ as either zero or one to select Dirichlet or Neumann boundary conditions. The diffusion constant is $k$, which is usually chosen as 1. The initial condition $u_0$ is a prescribed function. Throughout the study, we sue the heat equation as a preliminary benchmark since numerical methods for it are the easiest to implement.

Similarly, the wave equation is the classical example of a hyperbolic PDE used to model wave phenomena. Examples of physical systems modeled by the wave equation include acoustic waves, electromagnetic waves, and elastic waves.  We consider the wave equation 

\begin{gather}
    \begin{cases}
        u_{tt} = c^2\Delta u,& x\in \Omega, t> 0,\\
        u(0,x) = u_0(x),& x\in \Omega,\\
        u_t(0,x) = v_0(x),& x\in \Omega,\\
        \alpha_1 u + \alpha_2\pdv{u}{n} = 0,& x\in \partial \Omega, t>0,\\
        \beta_1 u + \beta_2 \pdv{u}{n} =0,& x\in \partial \Omega, t>0
    \end{cases}\label{eq:wave}
\end{gather}

\noindent
on the bounded domains on $\Omega = [-1,1], [0,1]$ and $[-1,1]^2$. The constants $\alpha_1,\alpha_2, \beta_1,$ and $\beta_2$ are chosen as before to select Dirichlet or Neumann boundary conditions. The wave speed is $c$, which is usually chosen as 1. The initial conditions $u_0$ and $v_0$ are again prescribed functions.

The Kuramoto Shivasinsky (KS) equation is used to model extended physical systems that have been driven far from equilibrium by intrinsic instabilities \citep[see][]{KSequ}. Examples of physical systems modeled by the KS equation include laminar flame fronts, reaction-diffusion systems, fluid films on inclines, and plasma dynamics. In the KS equation, the $u_{xx}$ term corresponds to large scale instabilities, the $u_{xxxx}$ corresponds to damping at small scales, and the nonlinear $uu_{x}$ term couples the small and large scales through energy transfer. Solutions to the KS equation become chaotic as the length of the domain is increased beyond $2\pi$. However, it is known that the attractor for the KS equation is compact with finite Hausdorff dimension, which motivates the use of manifold learning to create ROMs for the KS equation \citep[see][]{KSequ}. In this study, we use the equation
\begin{gather}
    \begin{cases}
    u_t = -uu_{x} - u_{xx} - u_{xxxx},& x\in (0, 32\pi), t> 0\\
    u(0,x) = u_0(x),& x\in(0,32\pi),
    \end{cases}\label{eq:KS}
\end{gather}
on the domain $[0, 32\pi]$ and consider periodic boundary conditions.

\subsection{Numerical Solutions}

To solve the full-order model for the heat equation, we use a Crank-Nicolson finite difference discretization \citep[see][]{heat-fd}

\begin{gather}
    \vb{A}\vb{u}^{n+1} = \vb{B}\vb{u}^n, \qqtext{where}\\
    \vb{A} = \mqty(2+2r& -r& 0&0& \dots& 0\\
    -r&2+2r&-r&0&\dots&0\\
    \ddots&\ddots&\ddots&\ddots&\ddots&\ddots\\
    0&\dots&0&-r&2+2r&-r\\
    0&\dots&0&0&-r&2+2r),\\ \vb{B} = \mqty(2-2r& r& 0&0& \dots& 0\\
    r&2-2r&r&0&\dots&0\\
    \ddots&\ddots&\ddots&\ddots&\ddots&\ddots\\
    0&\dots&0&r&2-2r&r\\
    0&\dots&0&0&r&2-2r), \\
    \qqtext{and} \quad r = \frac{\Delta t}{\Delta x^2}.
\end{gather}\label{eq:heat-fd}

For the wave and KS equations, we consider two finite-difference discretizations. First, we employ a summation by parts (SBP) finite-difference discretization \citep[see][]{ranocha2016summation}. This fourth-order discretization allows us to study more complicated initial conditions when we want to test the autoencoder. When we want to consider latent space dynamics, however, we switch to the following lower-order finite-difference scheme to ease the implementation (for Dirichlet boundary conditions) \citep[see][]{wave-fd}:

\begin{gather}
    (4/r^2 {\bf I} + {\bf K}){\bf u}^{n+1} =  2(4/r^2 {\bf I} - {\bf K}) {\bf u}^{n}  - (4/r^2 {\bf I} + {\bf K}){\bf u}^{n-1}
 , \nonumber\\\qqtext{where}
    {\bf K} = \left(\begin{array}{cccc}
2 & -1 & & \\
-1 & \ddots & \ddots & \\
 & \ddots & \ddots & -1 \\
 & & -1 & 2
\end{array}\right) ,
\qqtext{and}\quad r = \frac{c\Delta t}{\Delta x}. \label{eq:wave-fd}
\end{gather}

\begin{figure}[H]
     \centering
     \begin{subfigure}[b]{0.33\textwidth}
         \centering
         \includegraphics[width=\textwidth]{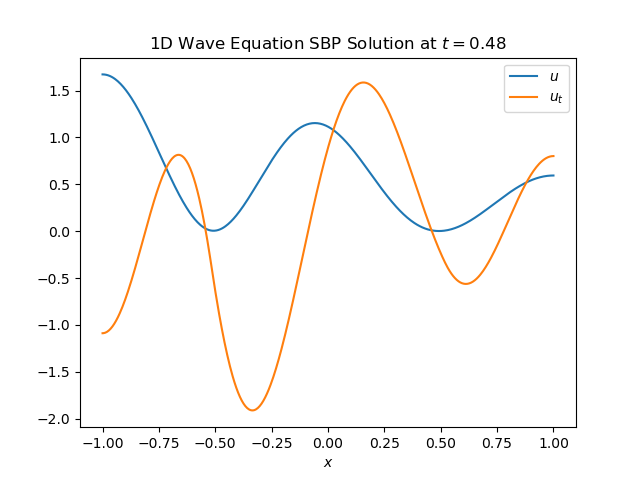}
         \caption{}
         \label{fig:wave_sol_1d}
     \end{subfigure}
     \hfill
     \begin{subfigure}[b]{0.32\textwidth}
         \centering
         \includegraphics[width=\textwidth]{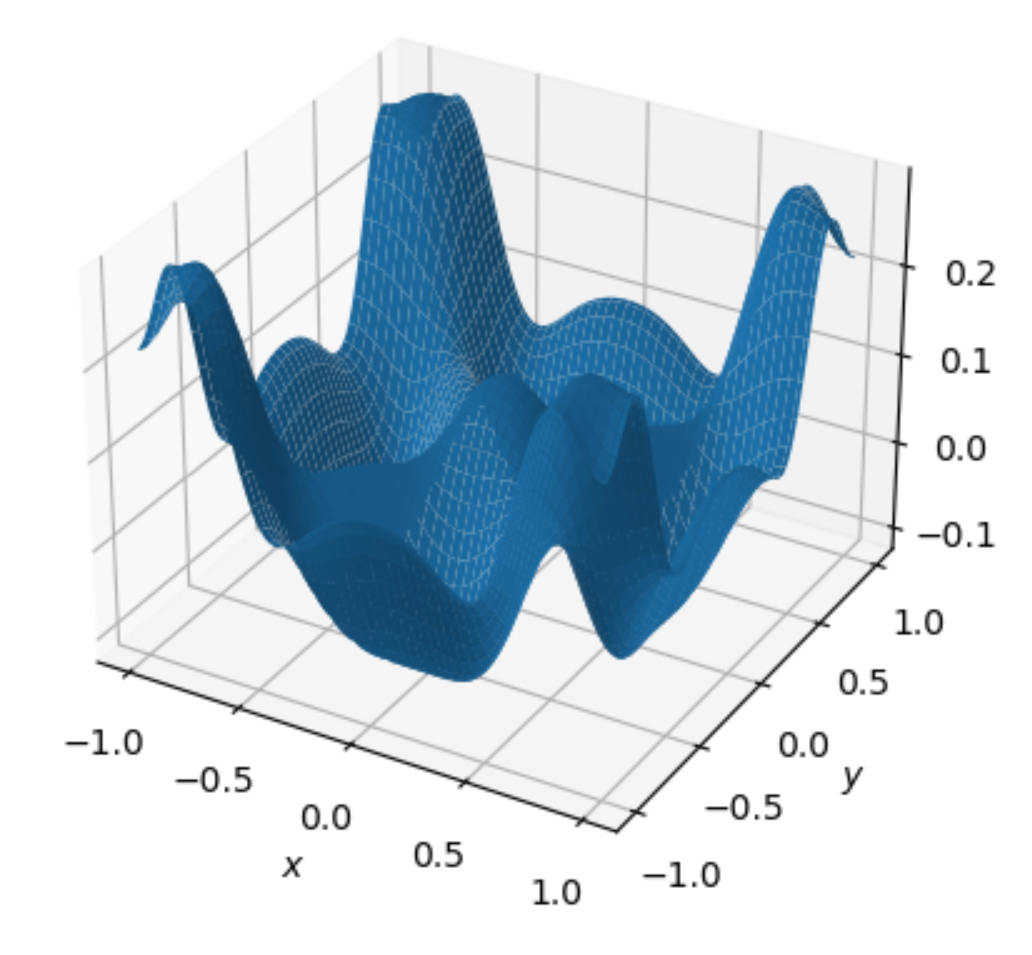}
         \caption{}
         \label{fig:wave_sol_2d}
     \end{subfigure}
     \hfill
     \begin{subfigure}[b]{0.32\textwidth}
         \centering
         \includegraphics[width=\textwidth]{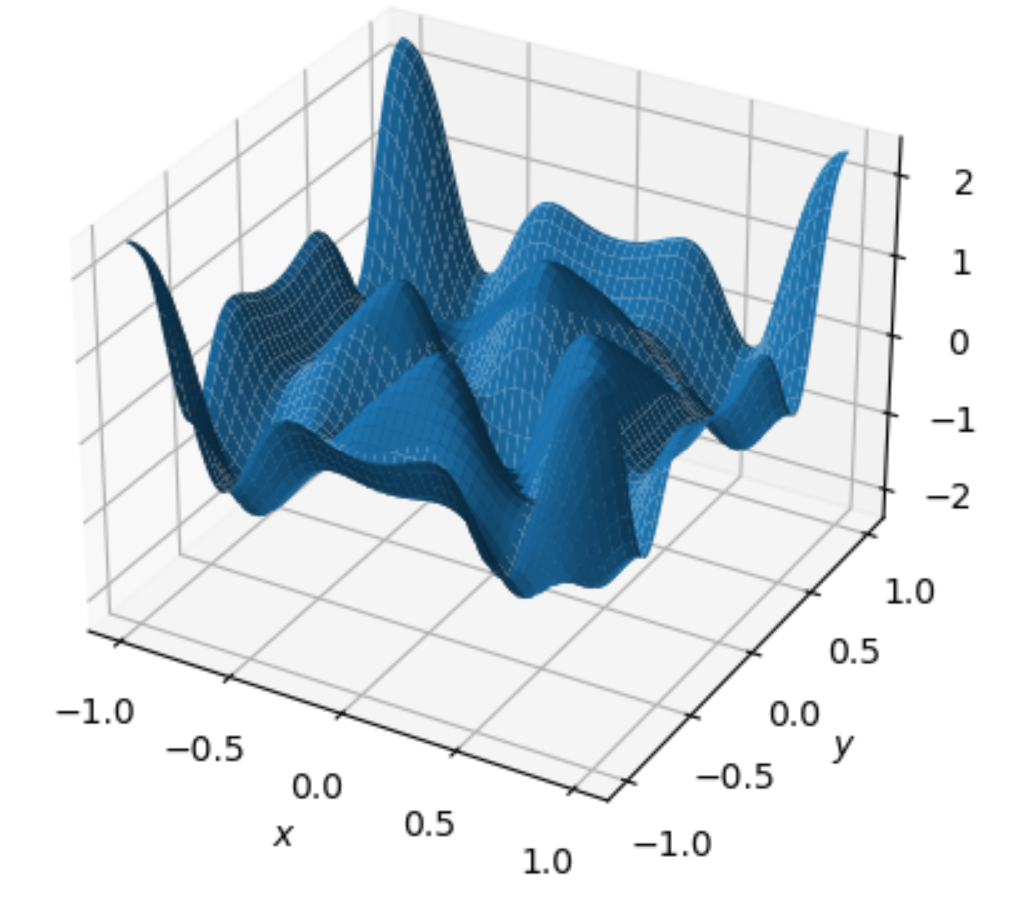}
         \caption{}
         \label{fig:wave_sol_2d_der}
     \end{subfigure}
        \caption{Numerical SBP solution to Wave Equation: (a) 1D, (b) 2D solution, and (c) 2D time derivative at $t=0.48$.}
        \label{fig:wave_solns}
\end{figure}

For the 2D wave equation, we again consider both a high-order SBP discretization and the following explicit finite-difference scheme  (again with Dirichlet boundary conditions) \citep[see][]{wave-fd}:

\begin{gather}
    \vb{u}^{n+1} = (2\vb{I} - \vb{K2D})\vb{u}^{n-1} - \vb{u}^{n-2},\\
    \qqtext{with} \quad \vb{K2D} = r^2\vb{I}\otimes \vb{K} + \vb{K} \otimes r^2\vb{I},
\end{gather}\label{eq:wave2d-fd}
and $r, \vb{K}$ given by \eqref{eq:wave-fd}. In Figure \ref{fig:wave_solns}, we show the numerical SBP solutions to the wave equation with homogeneous Neumann boundary condtions and $u_0(x) = e^{-x}\sin^2(\pi x)$ (1D), $u_0(x,y) = e^{-20(x^2+y^2)}$ (2D), with $v_0 \equiv 0$. For conciseness, the numerical results for the heat equation are presented in Section \ref{sec:dynamics}.

\begin{figure}[H]
         \centering
         \includegraphics[width=0.5\textwidth]{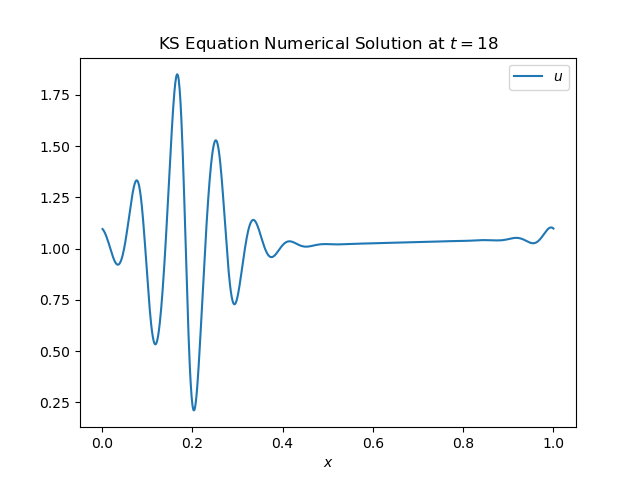}
         \caption{Fourier numerical solution of KS equation at $t=18$. }
         \label{fig:ks_soln}
\end{figure}

Similarly, for the Kuramoto-Shivasinsky equation, we also use two different discretizations. First, we consider a Fourier spatial discretization with 1024 modes and fourth-order Runge-Kutta time integration. Since the KS equation does not conserve any interesting quantities, we do not need to use an SBP discretization. Moreover, the high frequency oscillations are best captured by a high Fourier resolution. A plot of the numerical solution to the KS equation is given in Figure \ref{fig:ks_soln}. In this figure, we see behavior reminiscent of wavepackets --- there are both high and low frequency modes in the solution. The small oscillation at the right end of the domain is due to the periodic boundary conditions. 

We also consider the following finite-difference discretization for the KS equation: 
\begin{gather}
   \vb{A}^n\vb{u}^{n+1} = \vb{D}\vb{u}^{n-1},\\ \qqtext{where} \quad \vb{A}^n = \mqty(a + b_2^n + e& c_1^n&a& 0&0&0&\dots &0\\
   a + b_2^n& e&c_2^n&a&0&0&\dots&0\\
   a&b_2^n&e&c_3^n&a&0&\dots&0\\
   \ddots&\ddots&\ddots&\ddots&\ddots&\ddots&\ddots&\ddots\\
   0&\dots&0&a&b_{k-2}^n&e&c_{k-1}^n&a\\
   0&\dots&0&0&a&b_{k-1}^n&e&a+c_{k-1}^n\\
   0&\dots&0&0&0&a&b_k^n&a+e+c_{k-1}^n)\\
   \qqtext{and} \quad \vb{D} = \mqty(a' + b' + e'& a'& c&0&0&0&\dots &0\\
   a' + b'& e'&c'&a'&0&0&\dots&0\\
   a'&b'&e'&c'&a'&0&\dots&0\\
   \ddots&\ddots&\ddots&\ddots&\ddots&\ddots&\ddots&\ddots\\
   0&\dots&0&a'&b'&e'&c'&a'\\
   0&\dots&0&0&a'&b'&e'&a'+c'\\
   0&\dots&0&0&0&a'&b'&a'+e'+c')\\
   \qqtext{with}\quad  a = \frac{\Delta t}{2(\Delta x)^4},\quad b_i^n = \frac{\Delta t}{2(\Delta x)^2} - \frac{2\Delta t}{(\Delta x)^4} - \frac{\Delta t}{4\Delta x}u_{i-1}^n,\\
   c_i^n = \frac{\Delta t}{2(\Delta x)^2} - \frac{2\Delta t}{(\Delta x)^4} + \frac{\Delta t}{4\Delta x}u_{i+1}^n,\quad e = 1+\frac{2\Delta t}{(\Delta x)^4} - \frac{\Delta t}{(\Delta x)^2} \\
   a' = -a,\quad b' = c' =-\qty(\frac{\Delta t}{2(\Delta x)^2} - \frac{2\Delta t}{(\Delta x)^4}),\quad e' = 1-\frac{2\Delta t}{(\Delta x)^4} + \frac{\Delta t}{(\Delta x)^2}.
\end{gather}\label{eq:ks-fd}
This discretization is necessary since the Fourier discretization involves complex numbers, and our autoencoders are configured for real data \citep[see][]{KS-fd, KS-fd-2}. A plot of the Fourier numerical solution of the KS equation is given in Figure \ref{fig:ks_soln}.

\subsection{Data Matrices}
For the KS equation, which is first order in time, we have the data matrix
\begin{align}
    A = \mqty[\mqty{u_1(t_1)\\\vdots\\u_n(t_1)}&\mqty{u_1(t_2)\\\vdots\\u_n(t_2)} & \dots & \mqty{u_1(t_m)\\\vdots\\u_n(t_m)} ] \in \R^{n\times m}.\label{eq:data_matrix_ks}
\end{align}
of solution vectors sampled at different times. Here, we use an equal time spacing and sample 200 snapshots. 

For the wave equation, we have a second-order PDE in time. Thus, the state vector $\vb{a}(t_j)$ includes the spatially discretized solution $\vb{u}(t_j)$ as well as the time derivative $\vb{v}(t_j) := \vb{u_t}(t_j)$ concatenated into one column. We form the data matrix $A$ as
\begin{align}
    A = \mqty[\mqty{u_1(t_1)\\\vdots\\u_n(t_1)\\v_1(t_1)\\\vdots\\v_n(t_1)}&\mqty{u_1(t_2)\\\vdots\\u_n(t_2)\\v_1(t_2)\\\vdots\\v_n(t_2)} & \dots & \mqty{u_1(t_m)\\\vdots\\u_n(t_m)\\v_1(t_m)\\\vdots\\v_n(t_m)} ] \in \R^{2n\times m}.\label{eq:data_matrix_wave}
\end{align}

\begin{figure}[H]
     \centering
     \begin{subfigure}[b]{0.49\textwidth}
         \centering
         \includegraphics[width=\textwidth]{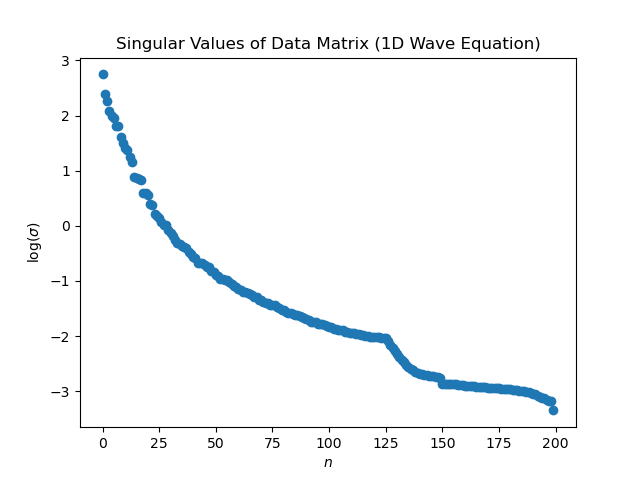}
         \caption{}
         \label{fig:wave_sigma}
     \end{subfigure}
     \hfill
     \begin{subfigure}[b]{0.49\textwidth}
         \centering
         \includegraphics[width=\textwidth]{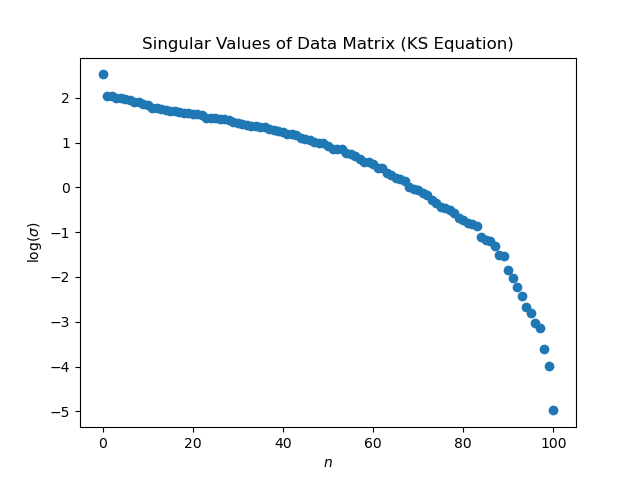}
         \caption{}
         \label{fig:ks_sigma}
     \end{subfigure}
        \caption{Singular values of data matrices for (a) 1D wave equation and (b) KS equation (log-scale).}
        \label{fig:svs}
\end{figure}

Once we have formed the data matrices, we can perform a preliminary analysis by plotting their singular values. Sharp decays or jumps in the singular values mean that linear projection ROM methods such as proper orthogonal decomposition (POD) can work well. A plot of the singular values for the KS and wave equations are given in Figure \ref{fig:svs}. From the figure, we can see that there are no significant gaps in the singular values for either system. This means that POD-type methods would not work well for creating ROMs. As a result, we turn to a \textit{nonlinear} projection method: manifold learning.

\section{Convolutional Autoencoders}
Convolutional autoencoders are types of artificial neural networks (ANNs) that allow us to ``learn'' lower-dimensional nonlinear manifolds that might describe our data. In our study, we found that due to the slow decay of the POD error, our data might not be well described by a lower-dimensional \textit{linear} subspace. As a result, we turn to convolutional autoencoders \citet{carlberg_ca}.

\begin{figure}[H]
    \centering
    \includegraphics[width=0.9\textwidth]{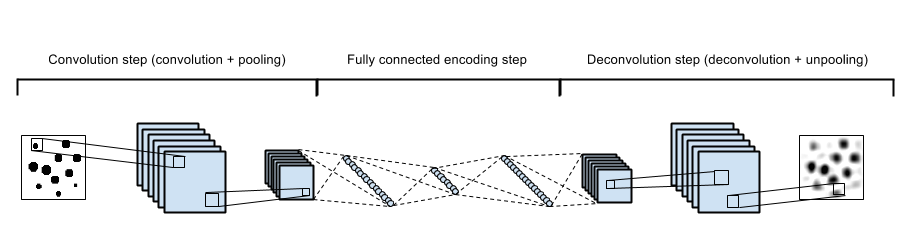}
    \caption{Covolutional Autoencoder structure \citep[see][]{conv_auto}.}
    \label{fig:conv_auto}
\end{figure}

Autoencoders function by mapping the high-dimensional input space $\R^n$ to itself in a way that factors through a lower-dimensional \textit{latent space} $\R^d$ as $\R^n \xrightarrow{f}\R^d\xrightarrow{g}\R^n$. The nonlinear functions $f$ and $g$ are known as the \textit{encoder} and \textit{decoder}, respectively, and they allow the input data to be approximated by passing to the latent space. The goal of the autoencoder is to construct the maps $f$ and $g$ such that the most salient features of the input data are captured in the latent space encoding of the data. The resultant approximations of the initial data form a differentiable nonlinear manifold known as the \textit{trial manifold}. Figure \ref{fig:conv_auto} shows a diagram of the convolutional autoencoder structure described above.

In this study, we build our convolutional autoencoders using the Keras and Tensorflow libraries for Python. These libraries provide a user interface with which to create customized neural networks by specifying the architecture and design parameters. The ANNs are described as graphs with vertices being the neurons and edges being the weights and activation functions. Most of the computations are carried out through C++ binaries for increased efficiency.

We build the autoencoder in two parts --- the encoder and decoder. Once we design the encoder architecture, we define the decoder architecture by reversing the order of the layers. For the encoder, we use $N_{\text{conv}}$ convolutional units followed by $N_{\text{dense}}$ fully connected layers. Each convolutional unit in the encoder consists of a convolutional layer followed by an AveragePooling layer, which reduces the dimension of the input by a factor of two. We use AveragePooling rather than MaxPooling to achieve a smoothing effect as we reduce dimension. Convolutional units in the decoder consist of convolutional layers followed by UpSampling layers, which ``invert'' the AveragePooling process to increase the dimension of the input by a factor of two. For all layers except the output layer of the decoder, we use a parametric rectified linear unit (PReLU) activation function, which has been shown to have advantages over the standard ReLU activation functions when working with signed data \citep[see][]{prelu}. For the output layer of the decoder, we use a linear activation function. We use a mean-squared error loss function with Adam as the optimizer \citep[see][]{adam}. 

\subsection{Generation of Training Data}

As with any ANN, it is important to have a large amount of training data to calibrate the model. One simple way to obtain this data would be to simply take snapshots of the state vector from the numerical solution of the full-order model. Essentially, this would amount to using the data matrix \eqref{eq:data_matrix_wave} or \eqref{eq:data_matrix_ks} as the training data. However, this goes against the fundamental purpose of creating reduced-order models; in order to train the reduced-order model, we should not need to evaluate the full-scale model hundreds or thousands of times. Instead, we employ a ``one size fits all'' approach. Note that the goal of the convolutional autoencoder is to learn the identity function $\text{id}_{\R^n}$. Thus, we only need to pass it sufficiently many vectors from $\R^n$. These vectors can in theory be \textit{completely unrelated} to the full-scale model to which we later apply the autoencoder. In this study, we use randomly generated smooth functions as our training data.

We devise two algorithms to generate these random smooth functions. The first method uses the Brownian bridge process (BBP). The BBP is a continuous stochastic process that creates a ``bridge'' between given points $p_1$ and $p_2$. Using a Python implementation of the BBP from \citet{brownian}, we obtain a Brownian bridge on the interval $[-1,1]$. This bridge can be thought of as a path from $x=-1$ to $x=1$ comprising $k$ steps with the $y$ value of the path varying randomly at each step. For each Brownian bridge used in the training data, we choose the number of steps $k$ randomly. 

Since we would like some regularity in our training data, we smoothen the Brownian bridges by applying a cubic spline interpolation on the BBP points. This interpolated function is evaluated at 500 points forming the spatial discretization of our numerical solution. Since our state vectors consist of two functions $u$ and $u_t$ concatenated together, we perform the BBP and interpolation twice to generate one element of the training dataset. An overview of the BBP/interpolation algorithm is given in Algorithm \ref{alg:bbp}.

While the Brownian Bridge process can be qualitatively interesting or ``clever'', it is computationally intensive. As a result, we turn to the far more efficient alternative: summing trigonometric functions. To generate a random smooth function $f_i$ using this method, we employ the formula
\begin{align*}
    f_i(x) = \sum_{j=1}^{N_i}A_i\sin(\omega_i x + \phi_i),
\end{align*}
where $N_i$ is sampled uniformly from $\qty{1, \dots, N_{\text{max}}}$, $A_i$ is sampled uniformly from $(0, A_{\text{max}}]$, $\omega_i$ is sampled uniformly from $(0, \omega_{\text{max}}]$, and $\phi_i$ is sampled uniformly from $[0, 2\pi)$. Since the trigonometric functions are smooth, we do not have to worry about interpolation. We used $N_{\text{max}} = 10, A_{\text{max}} = 5,$ and $\omega_{\text{max}}=10$ to generate our training data in the subsequent discussions. An overview of the Trigonometric algorithm is given in Algorithm \ref{alg:trig}.

\begin{algorithm}
\caption{BBP/Interpolation Algorithm}\label{alg:bbp}
\begin{algorithmic}[1]
\Procedure{BBP}{$N$, $\vb{x}$}\Comment{Generate $N$ random smooth functions over domain $\vu{x}$.}
\For{$i= 1:N$}
\State $n_x \gets$ \texttt{random-integer}(4, number of points in $\vb{x}$)\Comment{Minimum size of 4 so we can interpolate with cubic splines.}
\State $dt \gets 1/n_x$
\State $x_i \gets$ \texttt{random-normal}(0,1)
\State $B[0] \gets x_i$
\For{$n = 1:n_x-1$}
\State $t\gets n\cdot dt$
\State $x_i \gets \sqrt{dt}\cdot$ \texttt{random-normal}(0,1)
\State $B[n + 1] \gets B[n] \cdot (1 - dt / (1 - t)) + x_i$
\EndFor
\State $A[i] \gets$ interpolate $B$ over $(\vb{x}[0], \vb{x}[\text{end}])$
\EndFor
\State \textbf{return} $A$
\EndProcedure
\end{algorithmic}
\end{algorithm}

\begin{algorithm}
\caption{Trigonometric Algorithm}\label{alg:trig}
\begin{algorithmic}[1]
\Procedure{Trig}{$N$, $\vb{x}$, $N_{\text{max}}$, $\omega_{\text{max}}$, $A_{\text{max}}$}\Comment{Generate $N$ random smooth functions over domain $\vu{x}$ given parameters $N_{\text{max}}$, $\omega_{\text{max}}$, $A_{\text{max}}$.}
\For{$i= 1:N$}
\State $N_{i} \gets$ \texttt{random-integer}(1, $N_{\text{max}}$)
$A[i] \gets \vb{0}$
\For{$j=1:N_i$}
\State $\omega_{j} \gets$ \texttt{random-uniform}(0, $\omega_{\text{max}}$)
\State $A_{j} \gets$ \texttt{random-uniform}(0, $A_{\text{max}}$)
\State $\phi_{j} \gets$ \texttt{random-uniform}(0, $2\pi$)
\State $A[i] \gets A[i] + A_j\sin\qty(\omega_j\vb{x} + \phi_j)$
\EndFor
\EndFor
\State \textbf{return} $A$
\EndProcedure
\end{algorithmic}
\end{algorithm}

\begin{figure}[H]
     \centering
     \begin{subfigure}[b]{0.49\textwidth}
         \centering
         \includegraphics[width=\textwidth]{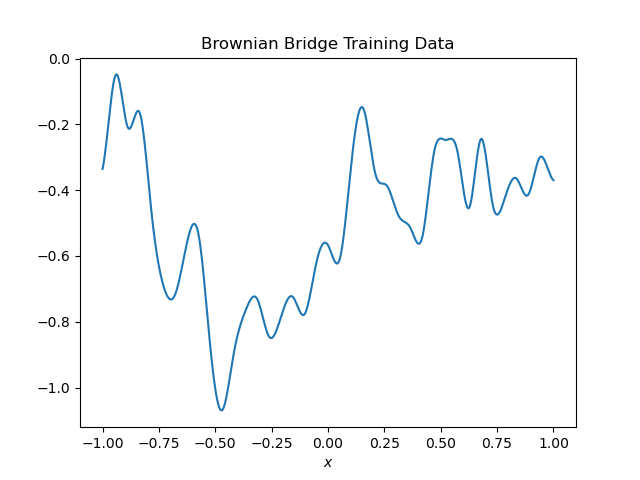}
         \caption{}
         \label{fig:bb_train}
     \end{subfigure}
     \hfill
     \begin{subfigure}[b]{0.49\textwidth}
         \centering
         \includegraphics[width=\textwidth]{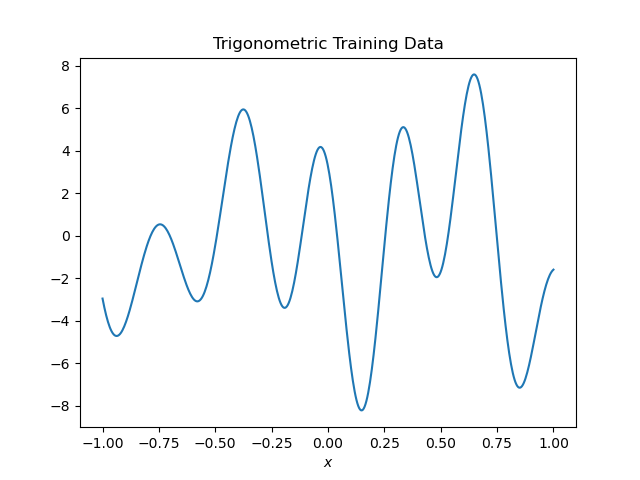}
         \caption{}
         \label{fig:trig_train}
     \end{subfigure}
        \caption{Training data samples: (a) Brownian Bridge/Interpolation and (b) Trigonometric.}
        \label{fig:train}
\end{figure}

For this study, we use a training dataset containing 60000 elements generated with one of the two algorithms described above. Figure \ref{fig:train} shows an example of one of the elements of the dataset generated with the BBP/interpolation method. We see that the stochastic nature of the BBP causes the training functions to have high frequency oscillations whereas our wave equation solutions do not. Overall, the generation of the training data took approximately 10 minutes. On the other hand, the trigonometric method takes less than 30 seconds to generate the same amount of data. A sample from the trigonometric training dataset is plotted in Figure \ref{fig:train}. 

\section{Results}
Using the trigonometric algorithm, we generated training data and proceeded to train the convolutional autoencoder using the Adam optimization routine in Keras (a variant of stochastic gradient descent with Nesterov momentum - \citet*{adam}). For training, we used a batch size of 128, a validation split of 20\%, and a learning rate of $\alpha = 0.01$. We trained the network until the validation loss reached $\sim 10^{-5}$, using the EarlyStopping callback in Keras to prevent overfitting. 

\begin{figure}[H]
     \centering
     \begin{subfigure}[b]{0.49\textwidth}
         \centering
         \includegraphics[width=\textwidth]{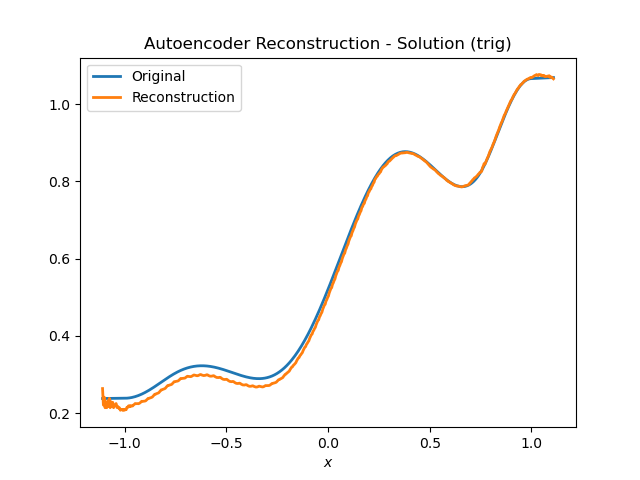}
         \caption{}
         \label{fig:wave_rec_s}
     \end{subfigure}
     \hfill
     \begin{subfigure}[b]{0.49\textwidth}
         \centering
         \includegraphics[width=\textwidth]{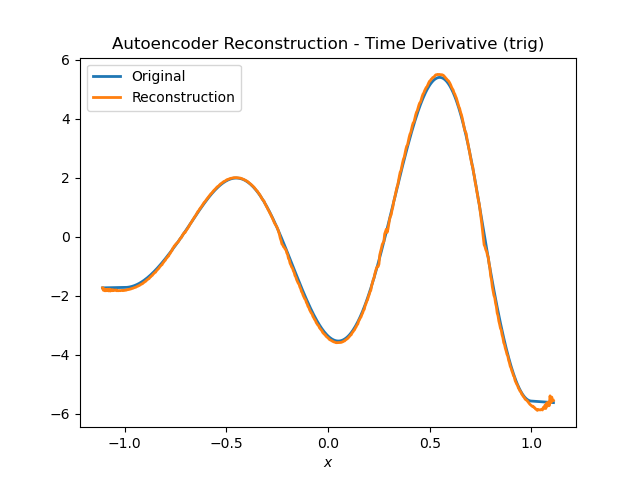}
         \caption{}
         \label{fig:wave_rec_d}
     \end{subfigure}
     \hfill
     \begin{subfigure}[b]{0.49\textwidth}
         \centering
         \includegraphics[width=\textwidth]{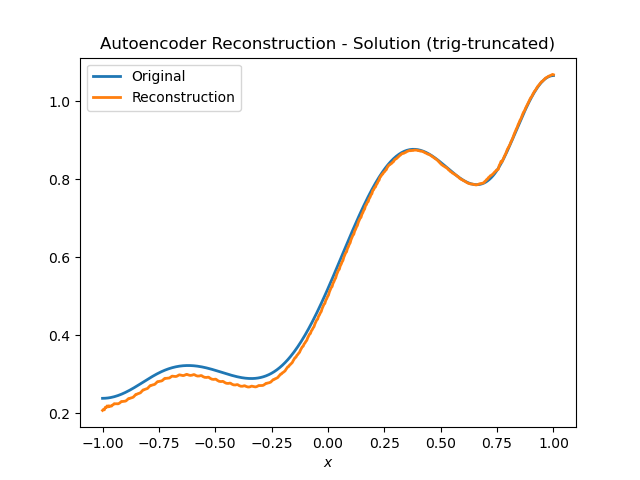}
         \caption{}
         \label{fig:wave_rec_s_t}
     \end{subfigure}
     \hfill
     \begin{subfigure}[b]{0.49\textwidth}
         \centering
         \includegraphics[width=\textwidth]{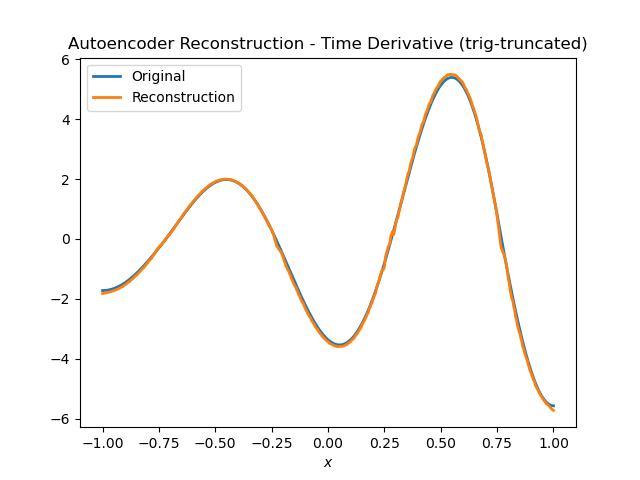}
         \caption{}
         \label{fig:wave_rec_d_t}
     \end{subfigure}
        \caption{Results of Convolutional Autoencoder --- reconstructions of 1D wave equation (a,c) solution and (b,d) time derivative at $t=0.48$. (a,b) are on the extended mesh and (c,d) are truncated to remove the boundary oscillations.}
        \label{fig:wave1d_recs}
\end{figure}

To view the results of the neural network training, we passed in snapshots of the numerical solution to the full-order models as test data. For a snapshot vector $\vb{u}$, the reconstructed vector is $\tilde{\vb{u}} = g(f(\vb{u}))$, where $f$ is the encoder and $g$ is the decoder. Plots of the original $\vb{u}$ versus reconstructed test data $\tilde{\vb{u}}$ are shown in Figure \ref{fig:wave1d_recs}. Here, we have dimension reduction by a factor of 32. From this figure, we see that the autoencoder does a very good job of reconstructing the test data. This is particularly interesting since the training data provided to the autoencoder was completely unrelated to the test data. We do notice that there are some Gibbs-phenomenon-type effects happening at the boundary points. We suspect that this may be due to the lack of boundary condition enforcement on the training data. 

Since the boundary irregularities quickly decay, we can avoid them using a mesh extension. We extend the functions to be defined on a slightly larger domain using tangent line approximations. Then, we truncate the reconstructions so that the boundary artifacts disappear. An example of this method is shown in \ref{fig:wave1d_recs}. In practice we tested mesh sizes of 100, 500, and 1000, finding that an extension of ~5\% of the original domain size works well.  

\begin{figure}[H]
     \centering
     \begin{subfigure}[b]{0.49\textwidth}
         \centering
         \includegraphics[width=\textwidth]{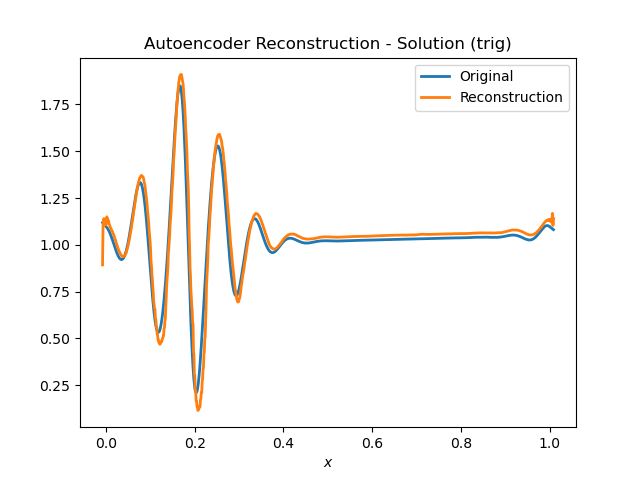}
         \caption{}
         \label{fig:ks_rec}
     \end{subfigure}
     \hfill
     \begin{subfigure}[b]{0.49\textwidth}
         \centering
         \includegraphics[width=\textwidth]{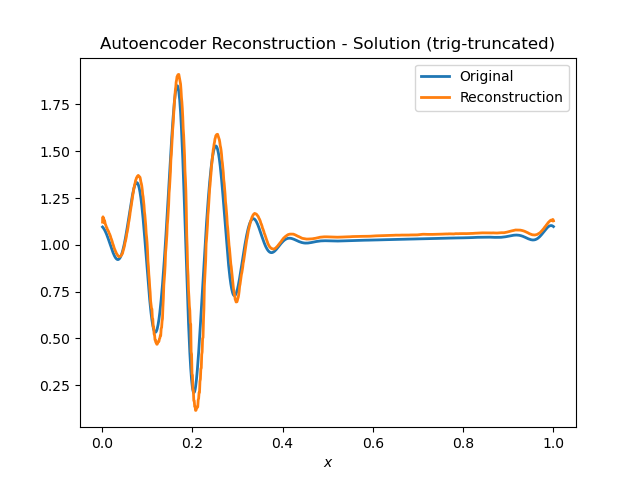}
         \caption{}
         \label{fig:ks_rec_t}
     \end{subfigure}
        \caption{Results of Convolutional Autoencoder --- reconstructions of KS equation solutions on (a) extended and (b) truncated mesh.}
        \label{fig:ks_recs}
\end{figure}

The results for the KS model are shown in Figure \ref{fig:ks_recs}. Here, we have dimension reduction by a factor of 16. We find that an autoencoder trained using the same method can work well for both the wave equation and the KS equation. The only parameter that needs to be changed in generating the training data is the size of the mesh. Here, we used mesh sizes of $500$ or $1000$ grid points for the 1D wave equation and $1024$ for the KS equation. This leads to the possibility of developing a ``bank'' of autoencoders trained on different mesh sizes (powers of 2, powers of 10, etc.). When needed, a pretrained autoencoder can be loaded from this bank and can be readily used. 

\begin{figure}[H]
     \centering
     \begin{subfigure}[b]{0.49\textwidth}
         \centering
         \includegraphics[width=\textwidth]{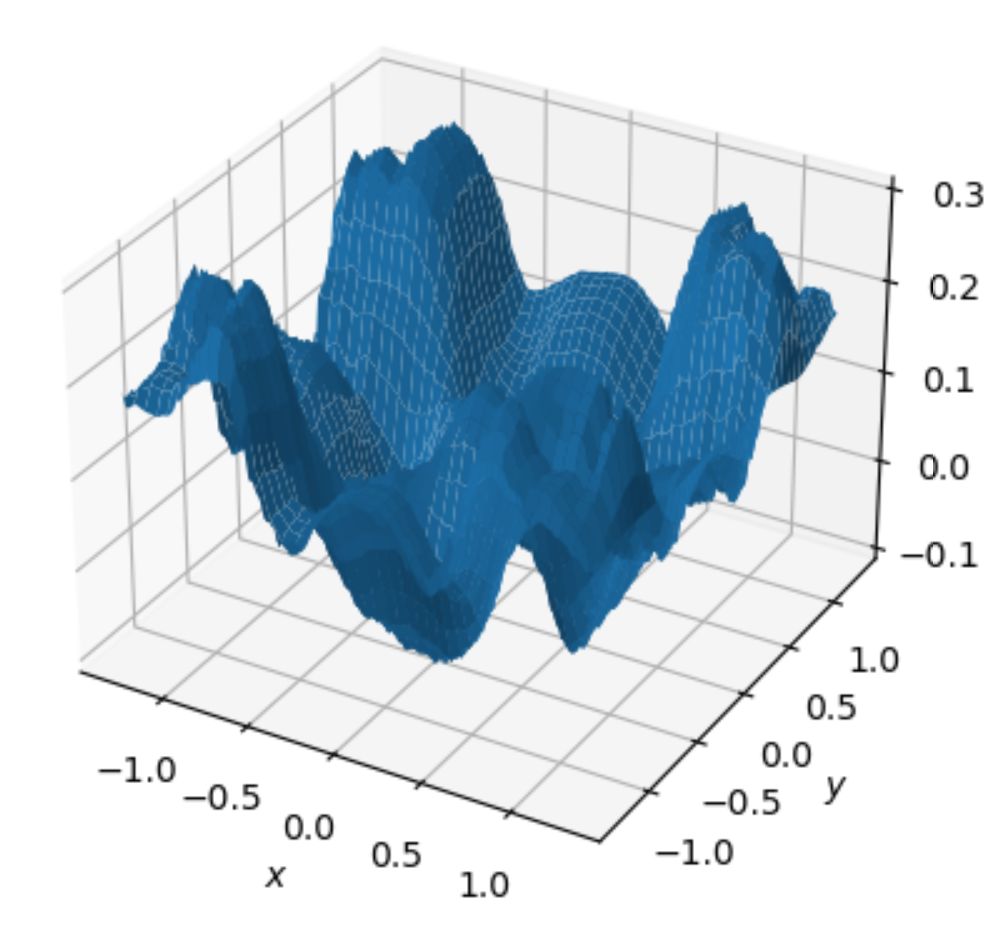}
         \caption{}
         \label{fig:wave2d_sol}
     \end{subfigure}
     \hfill
     \begin{subfigure}[b]{0.49\textwidth}
         \centering
         \includegraphics[width=\textwidth]{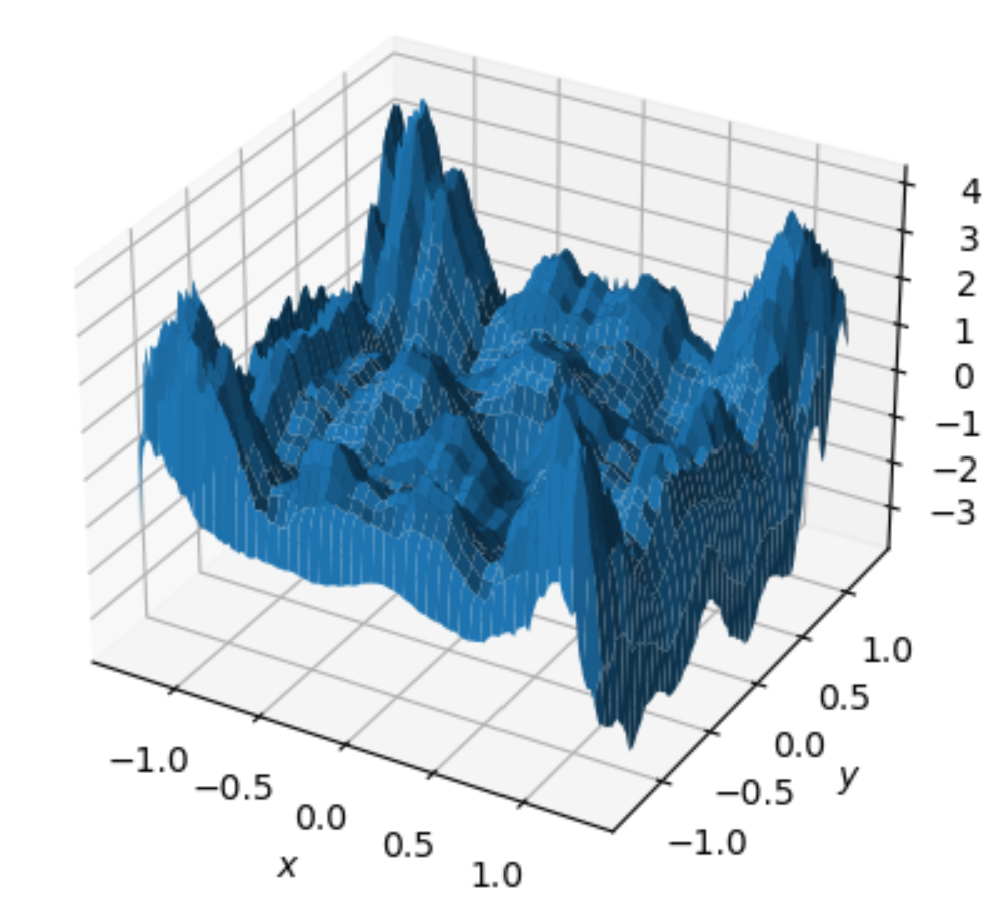}
         \caption{}
         \label{fig:wave2d_der}
     \end{subfigure}
     \hfill
     \begin{subfigure}[b]{0.49\textwidth}
         \centering
         \includegraphics[width=\textwidth]{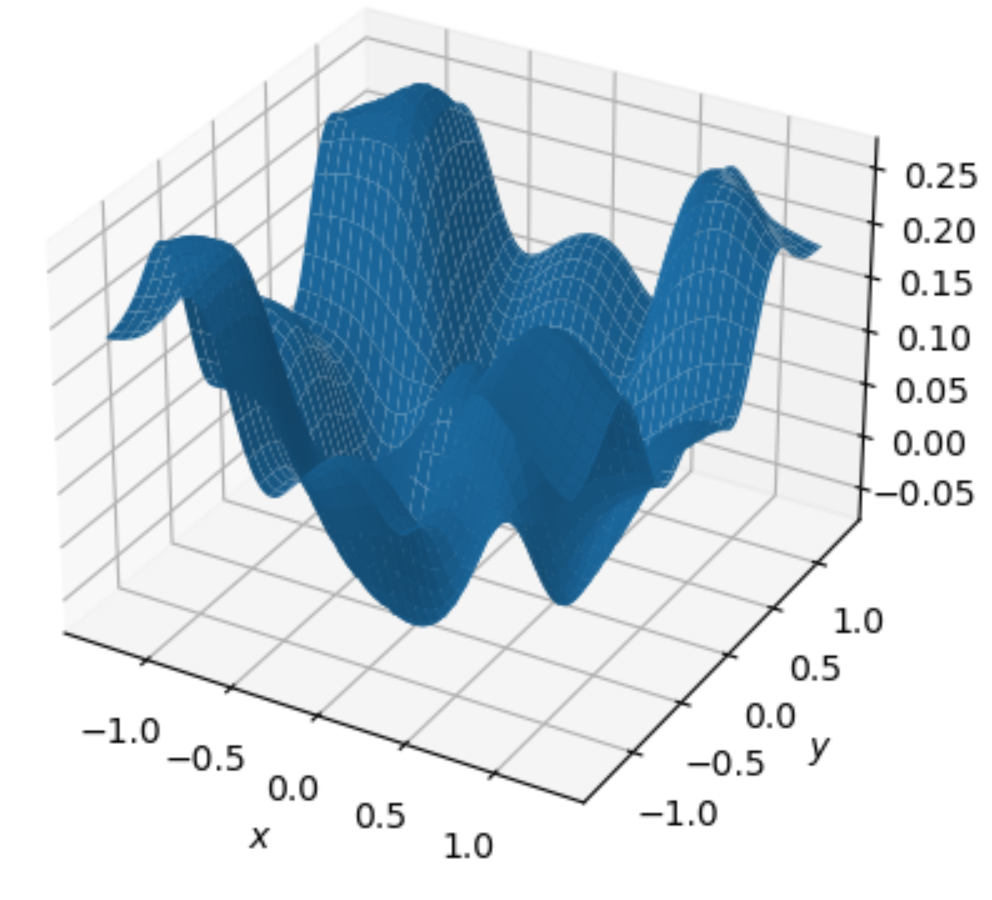}
         \caption{}
         \label{fig:wave2d_sol_blur}
     \end{subfigure}
     \hfill
     \begin{subfigure}[b]{0.49\textwidth}
         \centering
         \includegraphics[width=\textwidth]{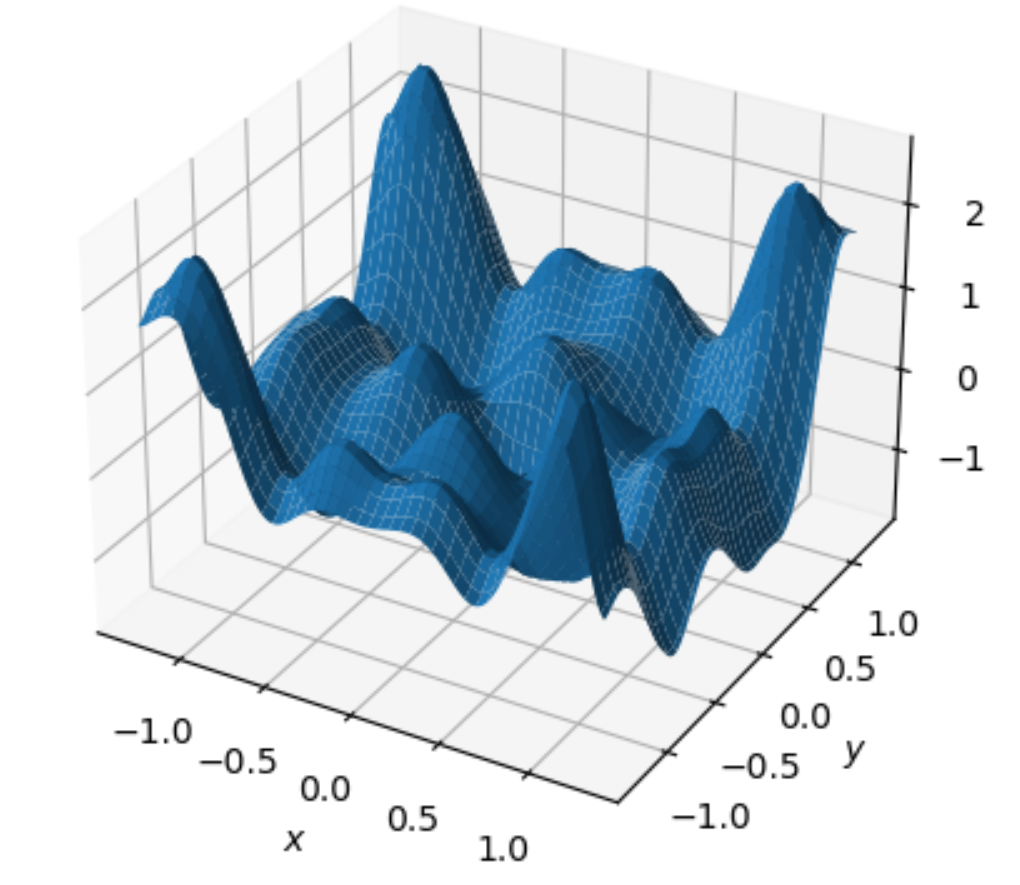}
         \caption{}
         \label{fig:wave2d_der_blur}
     \end{subfigure}
        \caption{Results of convolutional autoencoder --- reconstructions of 2D wave equation (a,c) solution and (b,d) time derivative at $t=0.48$. (a,b) are without Gaussian blur and (c,d) have Gaussian blur applied. All plots are truncated to remove the boundary oscillations.}
        \label{fig:wave2d_recs}
\end{figure}

We also tested the autoencoder on the 2D wave equation. To generate training data for this case, we simply tensored the 1D training data --- that is, we took $f(x,y) = g(x)h(y)$, where $g$ and $h$ are random 1D functions generated with the methods described previously. When training the autoencoder on the 2D data for large grid sizes, we can run into issues of high memory usage. There are two ways to avoid this problem. The first is to dynamically load the training data in chunks using \textit{generators}. However, when working with data that is not in the form of images, this method is difficult to implement directly in Keras. We found that an alternative solution is to reduce the size of training data. In fact, we saw that using fewer training data increased the training time, but did not significantly impact the final accuracy of the autoencoder, as shown in table Table \ref{tab:train}. The results of the 2D case are plotted in Figure \ref{fig:wave2d_recs}. Again, we circumvent the boundary oscillations by extending the functions using a tangent plane approximation. In the 2D reconstuctions, we also find some oscillatory behavior in the interior of the domain. To smoothen out the reconstructions, we can apply a Gaussian filter as needed \citep[see][]{gaussian-filter}. Moreover, we can add the Gaussian filter to the final layer of the decoder to achieve the same effect. It is also worth noting that the core architecture of the autoencoder remains the same for 1D and 2D problems. The only structural changes necessary are the changes in input and output dimensions. In addition, the autoencoder achieves greater dimension reduction on 2D problems since the reduction factor is squared. Here, we have a reduction factor of $8^2=64$.

\begin{table}[H]
\centering
\resizebox{\textwidth}{!}{%
\begin{tabular}{@{}lll@{}}
\toprule
\multicolumn{1}{c}{Size of Training Data} & \multicolumn{1}{c}{Average Number of Training Epochs} & \multicolumn{1}{c}{Final Loss Function Value} \\ \midrule
500    & 17 & $\sim$$10^{-4}$ \\
1,000  & 12 & $\sim$$10^{-4}$ \\
10,000 & 6  & $\sim$$10^{-4}$ \\
20,000 & 4  & $\sim$$10^{-4}$ \\ \bottomrule
\end{tabular}%
}
\caption{Training performance versus size of training data}
\label{tab:train}
\end{table}

\subsection{Latent Space Dynamics: Manifold LSPG}\label{sec:dynamics}

Next, we consider the manifold least-squares Petrov-Galerkin projection method outlined in \citet{carlberg_ca} to carry out the latent space dynamics. For simplicity, we first consider it on the 1D heat equation

The idea of the manifold LSPG method is to minimize the norm of the time discrete residual over the latent space. That is, we obtain the encoded solution vector at the $n^{\text{th}}$ timestep by solving
\begin{align}
    \vu{u}^n = \min_{\vb{\xi}\in \R^d} \norm{\mathbf{r}^n\qty(\mathbf{x}_{\text{ref}} + g(\xi))}_2^2.\label{eq:opt}
\end{align}
Here, $\mathbf{r}^n$ is the time-discrete residual for the ODE describing the evolution of the spatially discretized full-order model. Since we use a Crank-Nicolson time discretization for the full-order model \eqref{eq:heat-fd}, we obtain
\begin{align*}
  \mathbf{r}^n(\vb{y}) = \vb{A}\vb{y} - \vb{B}\vb{u}^{n-1},
\end{align*}
where $\vb{u}^{n=1}$ is the solution at the previous timestep and $f$ is the RHS of the ODE from the full-order model. In \eqref{eq:opt}, $g$ is the decoder (coming from the autoencoder) and
\begin{align*}
    \mathbf{x}_{\text{ref}} = \vb{u}_0 - \texttt{autoencoder}(\vb{u}_0)
\end{align*}
is the reference value that sets the ``origin'' for the coordinates on the lower-dimensional manifold. We minimize over vectors in $\R^d$, where $d$ is the dimension of the latent space.

To carry out the optimization, we employ a Gauss-Newton algorithm as described in Section 3.4 of \citet{carlberg_ca}. For efficiency, we use the built-in methods from $\texttt{scipy.optimize}$. To compute the Jacobian of the decoder, we use the automatic differentiation tools in Tensorflow. 

We find that the optimization works very well in obtaining the latent space dynamics (Figure \ref{fig:lspg}). To test the robustness of the method, we then used the same autoencoder to compute an LSPG projection for the 1D wave equation \eqref{eq:wave} with speec $c=10$, initial condition $u_0(x) = 5x(1-x)$ and $v_0(x) = 3\sin(\pi x)$ and homogeneous Dirichlet boundary conditions. To compute the LSPG, we again used the Crank-Nicolson finite-difference scheme in \eqref{eq:wave-fd}.

This gives the residual
\begin{align*}
    \mathbf{r}^n(\vb{y}) = (4/r^2 {\bf I} + {\bf K}){\bf y}-2(4/r^2 {\bf I} - {\bf K}) {\bf u}^{n-1}  + (4/r^2 {\bf I} + {\bf K}){\bf u}^{n-2}.
\end{align*}

Next, we use the same autoencoder to construct the LSPG ROM for the KS equation \eqref{eq:KS} with initial condition $u_0(x) = \cos(\pi x/16)$ and boundary conditions $u_x(t,0) = 0$ (all others free) on the domain $[0,32\pi]$. Since our previous discretization for the KS equation was based on a Fourier decomposition, we cannot directly implement the LSPG algorithm for the KS equation (the Fourier coefficients are complex-valued). As a result, we use a finite-difference scheme for the LSPG ROM \eqref{eq:ks-fd}, which gives the residual
\begin{align*}
    \mathbf{r}^n(\vb{y}) = \vb{A}\vb{y} - \vb{D}\vb{u}^{n-1}.
\end{align*}

Finally, for the 2D wave equation discretized with the finite-difference scheme \eqref{eq:wave2d-fd}, giving the residual

\begin{align*}
    \mathbf{r}^n(\vb{y}) = \vb{y} - (2\vb{I} - \vb{K2D})\vb{u}^{n-1} + \vb{u}^{n-2}.
\end{align*}

The results of the LSPG projections are given in Figure \ref{fig:lspg}. From this figure, we see that an LSPG projection with the \textit{same} autoencoder gives excellent results for three completely different dynamical systems. This again motivates the creation of an ``autoencoder bank'' that contains trained autoencoders that can be readily deployed on several real-world problems. 

\begin{figure}[H]
     \centering
     \begin{subfigure}[b]{0.49\textwidth}
         \centering
         \includegraphics[width=\textwidth]{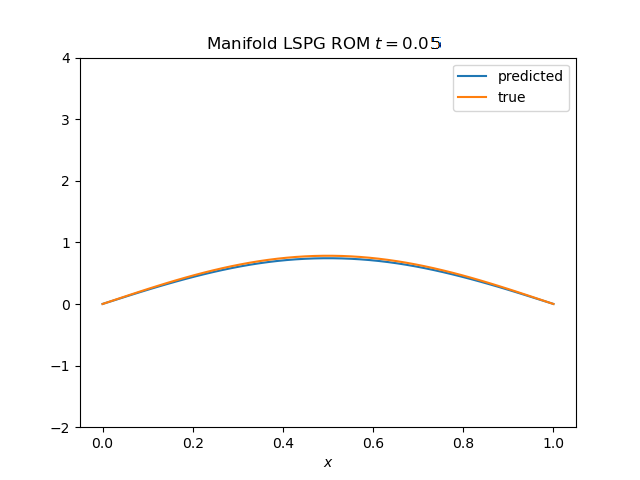}
         \caption{}
         \label{fig:heat_rom}
     \end{subfigure}
     \hfill
     \begin{subfigure}[b]{0.49\textwidth}
         \centering
         \includegraphics[width=\textwidth]{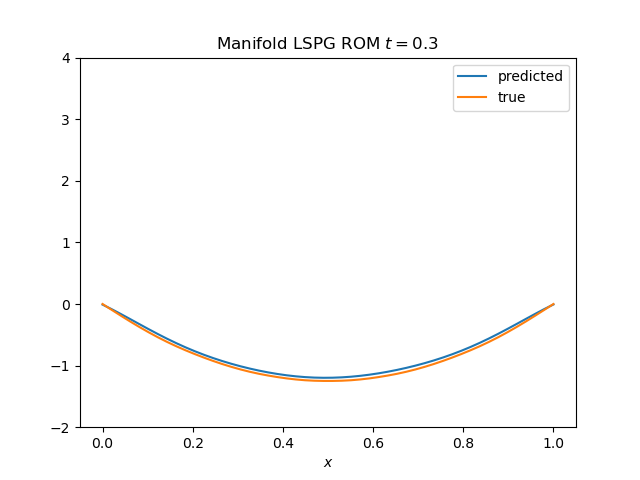}
         \caption{}
         \label{fig:wave_rom}
     \end{subfigure}
\end{figure}
\begin{figure}[H]\ContinuedFloat
\centering
     \begin{subfigure}[b]{0.49\textwidth}
         \centering
         \includegraphics[width=\textwidth]{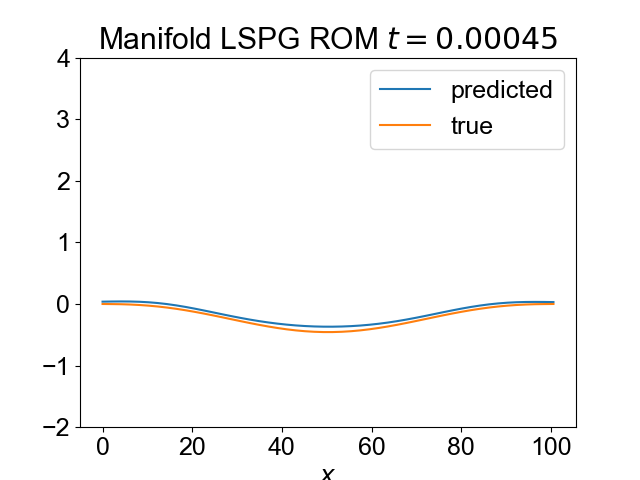}
         \caption{}
         \label{fig:KS_rom}
     \end{subfigure}
     \hfill
     \begin{subfigure}[b]{0.49\textwidth}
         \centering
         \includegraphics[width=\textwidth]{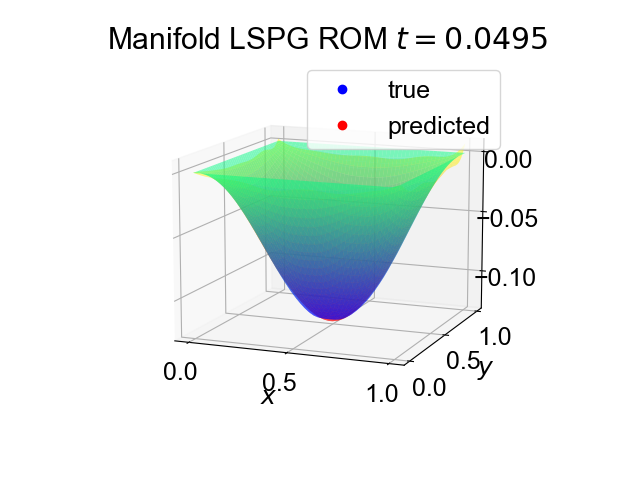}
         \caption{}
         \label{fig:wave2d_rom}
     \end{subfigure}
\end{figure}
\begin{figure}[H]\ContinuedFloat
\centering
     \begin{subfigure}[b]{0.49\textwidth}
         \centering
         \includegraphics[width=\textwidth]{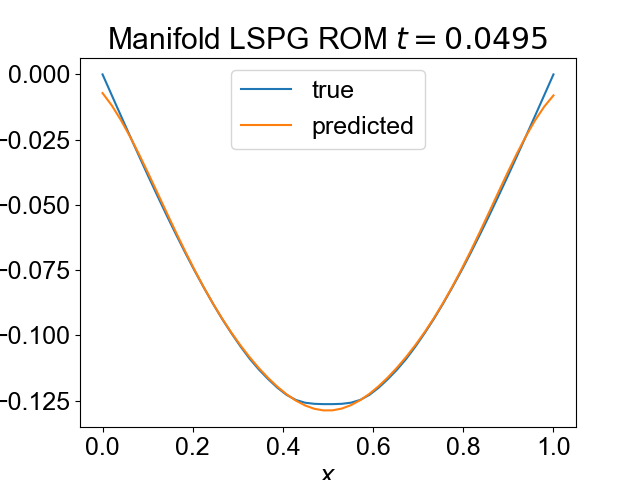}
         \caption{}
         \label{fig:wave2d_slice_rom}
     \end{subfigure}
        \caption{Results of manifold LSPG ROM --- ROM and actual solutions for (a) heat equation after 100 timesteps, (b) 1D wave equation after 600 timesteps, (c) KS equation after 50 timesteps, (d) 2D wave equation after 100 timesteps, and (e) 2D wave equation after 100 timesteps (slice at maximum amplitude).}
        \label{fig:lspg}
\end{figure}

For both of the cases in Figure \ref{fig:lspg}, we started with an extended mesh size of 64 and reduced it to 8 in the latent space. The largest computational cost in performing the LSPG projection is the calculation of the Jacobian of the decoder. However, the necessary backpropagation graphs can be precomputed in Tensorflow, so this is just a one-time overhead cost. In addition, for larger problems (such as the 2D wave equation), the decoder Jacobian may be too large to compute and store with the automatic differentiation method. In such cases using the ``3-point'' Jacobian option in \texttt{scipy.optimize.least\_squares} will work instead. The 3-point Jacobian gives comparable results to the automatic differentiation method for all the model problems described above. There is a tradeoff between time and memory: the automatic differentiation can be faster (if precomputed and stored), but uses much more memory. On the other hand, the 3-point Jacobian is slower (for larger problems) but does not require as much memory. All other matrix solves are done on the low-dimensional latent space. Here, we reduced from 64 to 8 dimensions for the 1D models and from $72\times 72$ dimensions to $9\times 9$ dimensions for the 2D wave equation.

\section{Conclusions and Future Work}
In this study, we considered the creation of reduced-order models based on nonlinear projection methods and manifold learning. We developed and trained a convolutional autoencoder to project the numerical solution to the wave equation onto a lower-dimensional nonlinear manifold. To generate the training data for the autoencoder, we used a combination of the Brownian Bridge Process and interpolation as well as a trigonometric process that was completely independent of the testing data. Regardless, we found that the autoencoder reconstructs the testing data extremely well. 

We also found that the manifold LSPG projection method works well to simulate the latent-space dynamics and that a single autoencoder can be used to create ROMs for completely different dynamical systems. Incorporating this projection technique with our method of creating convolutional autoencoders prompts the creation of a ``bank'' of autoencoder ROMs pre-trained on various mesh sizes that can be readily available for use in real-world applications.

Future extensions to these types of autoencoders would involve augmenting the training data generation process and training parameters to reduce the reconstruction error at the boundaries. For example, one could consider regularization during training and enforcement of boundary conditions on the training data. Additional extensions can use these autoencoders to define projection-based reduced-order model that conserves energy. In \citet{carlberg_cons}, the authors project a Finite Volume discretization via a convolutional autoencoder. For the wave equation example, we can consider projecting the SBP discretization onto a lower dimensional manifold. 

It is also interesting to consider the case where we are presented with data without a specific definition of an associated conserved quantity. Using autoencoders and symmetries in the data, \citet{noether} has developed a method to learn conservation laws from data. In addition, \citet{hnn} have used Hamiltonian Neural Networks to define reduced-order models based on learned conservation laws. These methods can potentially be used in conjunction with the convolutional autoencoders presented here to extract new physics from data and create fast, accurate reduced-order models.
\nocite{*}

\acks{We would like to acknowledge support for this project
from the National Science Foundation grants DMS-1745654 and DMS-1906446. The research of RCS was also supported in part by the Air Force Office of Scientific Research (AFOSR) through the grant FA9550-18-1-0457.}\\

\noindent Source code for this work can be obtained by contacting Sreeram Venkat (srvenkat@ncsu.edu). We are working on making the source code publically available on github.


\newpage

\appendix







\bibliography{refs.bib}

\begin{thebibliography}{25}
\providecommand{\natexlab}[1]{#1}
\providecommand{\url}[1]{\texttt{#1}}
\expandafter\ifx\csname urlstyle\endcsname\relax
  \providecommand{\doi}[1]{doi: #1}\else
  \providecommand{\doi}{doi: \begingroup \urlstyle{rm}\Url}\fi

\bibitem[Akrivis(1992)]{KS-fd-2}
Georgios~D Akrivis.
\newblock Finite difference discretization of the kuramoto-sivashinsky
  equation.
\newblock \emph{Numerische Mathematik}, 63\penalty0 (1):\penalty0 1--11, 1992.

\bibitem[Antoulas et~al.(2020)Antoulas, Beattie, and
  G{\"u}{\u{g}}ercin]{interpolatory}
Athanasios~Constantinos Antoulas, Christopher~Andrew Beattie, and Serkan
  G{\"u}{\u{g}}ercin.
\newblock \emph{Interpolatory methods for model reduction}.
\newblock SIAM, 2020.

\bibitem[Benner et~al.(2015)Benner, Gugercin, and Willcox]{projection}
Peter Benner, Serkan Gugercin, and Karen Willcox.
\newblock A survey of projection-based model reduction methods for parametric
  dynamical systems.
\newblock \emph{SIAM review}, 57\penalty0 (4):\penalty0 483--531, 2015.

\bibitem[Bui-Thanh et~al.(2008)Bui-Thanh, Willcox, and
  Ghattas]{high-dim-model-reduction}
Tan Bui-Thanh, Karen Willcox, and Omar Ghattas.
\newblock Model reduction for large-scale systems with high-dimensional
  parametric input space.
\newblock \emph{SIAM Journal on Scientific Computing}, 30\penalty0
  (6):\penalty0 3270--3288, 2008.

\bibitem[Deng and Cahill(1993)]{gaussian-filter}
Guang Deng and LW~Cahill.
\newblock An adaptive gaussian filter for noise reduction and edge detection.
\newblock In \emph{1993 IEEE conference record nuclear science symposium and
  medical imaging conference}, pages 1615--1619. IEEE, 1993.

\bibitem[Dong(2006)]{wave-fd}
Shuonan Dong.
\newblock Finite difference methods for the hyperbolic wave partial
  differential equations.
\newblock \emph{Computational Science and Engineering II}, 2006.

\bibitem[Galbally et~al.(2010)Galbally, Fidkowski, Willcox, and
  Ghattas]{galbally2010non}
David Galbally, Krzysztof Fidkowski, Karen Willcox, and Omar Ghattas.
\newblock Non-linear model reduction for uncertainty quantification in
  large-scale inverse problems.
\newblock \emph{International journal for numerical methods in engineering},
  81\penalty0 (12):\penalty0 1581--1608, 2010.

\bibitem[Gao et~al.(2018)Gao, Cai, and Liu]{KS-fd}
Ping Gao, Chengjian Cai, and Xiaoyi Liu.
\newblock Numerical simulation of stochastic kuramoto-sivashinsky equation.
\newblock \emph{Journal of Applied Mathematics and Physics}, 6\penalty0
  (11):\penalty0 2363--2369, 2018.

\bibitem[Gramacy(2020)]{gaussian-processes}
Robert~B Gramacy.
\newblock \emph{Surrogates: Gaussian process modeling, design, and optimization
  for the applied sciences}.
\newblock Chapman and Hall/CRC, 2020.

\bibitem[Greydanus et~al.(2019)Greydanus, Dzamba, and Yosinski]{hnn}
Samuel Greydanus, Misko Dzamba, and Jason Yosinski.
\newblock Hamiltonian neural networks.
\newblock In \emph{Advances in Neural Information Processing Systems}, pages
  15379--15389, 2019.

\bibitem[Jones(2015)]{conv_auto}
Swarbrick Jones.
\newblock Convolutional autoencoder, 2015.
\newblock URL \url{https://images.app.goo.gl/HXcwZPe2YQVFrUup7}.

\bibitem[Kunisch and Volkwein(1998)]{POD}
K~Kunisch and S~Volkwein.
\newblock Control of burgers’ equation by a reduced order approach using
  proper orthogonal decomposition, optimierung und kontrolle bericht 138
  (1998), universitat graz, austria.
\newblock \emph{J. Opt. Theory Applic., to appear}, 1998.

\bibitem[Kunisch and Volkwein(2010)]{kunisch2010optimal}
Karl Kunisch and Stefan Volkwein.
\newblock Optimal snapshot location for computing pod basis functions.
\newblock \emph{ESAIM: Mathematical Modelling and Numerical
  Analysis-Mod{\'e}lisation Math{\'e}matique et Analyse Num{\'e}rique},
  44\penalty0 (3):\penalty0 509--529, 2010.

\bibitem[Lee and Carlberg(2019)]{carlberg_cons}
Kookjin Lee and Kevin Carlberg.
\newblock Deep conservation: A latent dynamics model for exact satisfaction of
  physical conservation laws.
\newblock \emph{arXiv preprint arXiv:1909.09754}, 2019.

\bibitem[Lee and Carlberg(2020)]{carlberg_ca}
Kookjin Lee and Kevin~T Carlberg.
\newblock Model reduction of dynamical systems on nonlinear manifolds using
  deep convolutional autoencoders.
\newblock \emph{Journal of Computational Physics}, 404:\penalty0 108973, 2020.

\bibitem[Mototake(2019)]{noether}
Yoh-ichi Mototake.
\newblock Interpretable conservation law estimation by deriving the symmetries
  of dynamics from trained deep neural networks.
\newblock \emph{arXiv preprint arXiv:2001.00111}, 2019.

\bibitem[Oono(2016)]{brownian}
Kenta Oono.
\newblock Brownian bridge sample code, 2016.
\newblock URL \url{https://gist.github.com/delta2323/6bb572d9473f3b523e6e}.

\bibitem[Quarteroni et~al.(2015)Quarteroni, Manzoni, and
  Negri]{quarteroni2015reduced}
Alfio Quarteroni, Andrea Manzoni, and Federico Negri.
\newblock \emph{Reduced basis methods for partial differential equations: an
  introduction}, volume~92.
\newblock Springer, 2015.

\bibitem[Ranocha et~al.(2016)Ranocha, {\"O}ffner, and
  Sonar]{ranocha2016summation}
Hendrik Ranocha, Philipp {\"O}ffner, and Thomas Sonar.
\newblock Summation-by-parts operators for correction procedure via
  reconstruction.
\newblock \emph{Journal of Computational Physics}, 311:\penalty0 299--328,
  2016.

\bibitem[Rasmussen(2003)]{gaussian-processes-ML}
Carl~Edward Rasmussen.
\newblock Gaussian processes in machine learning.
\newblock In \emph{Summer school on machine learning}, pages 63--71. Springer,
  2003.

\bibitem[Rehmann and Janisch(2020)]{KSequ}
Ulf Rehmann and Maximillian Janisch.
\newblock Kuramoto-sivashinsky equation, 2020.
\newblock URL
  \url{https://encyclopediaofmath.org/wiki/Kuramoto-Sivashinsky_equation}.

\bibitem[Smith et~al.(1985)Smith, Smith, and Smith]{heat-fd}
Gordon~D Smith, Gordon~D Smith, and Gordon Dennis~Smith Smith.
\newblock \emph{Numerical solution of partial differential equations: finite
  difference methods}.
\newblock Oxford university press, 1985.

\bibitem[Zhang et~al.(2018)Zhang, Pan, Sun, and Tang]{prelu}
Yu-Dong Zhang, Chichun Pan, Junding Sun, and Chaosheng Tang.
\newblock Multiple sclerosis identification by convolutional neural network
  with dropout and parametric relu.
\newblock \emph{Journal of computational science}, 28:\penalty0 1--10, 2018.

\bibitem[Zhang(2018)]{adam}
Zijun Zhang.
\newblock Improved adam optimizer for deep neural networks.
\newblock In \emph{2018 IEEE/ACM 26th International Symposium on Quality of
  Service (IWQoS)}, pages 1--2. IEEE, 2018.

\bibitem[Zhao et~al.(2017)Zhao, Yan, Wang, and Mao]{lstm}
Rui Zhao, Ruqiang Yan, Jinjiang Wang, and Kezhi Mao.
\newblock Learning to monitor machine health with convolutional bi-directional
  lstm networks.
\newblock \emph{Sensors}, 17\penalty0 (2):\penalty0 273, 2017.

\end{thebibliography}

\end{document}